\documentstyle[12pt]{article}
\newtheorem{defi}{Definition}

\def\today{\ifcase\month\or
January\or February\or March\or April\or May\or June\or July\or August\or September\or October\or November\or December\fi
\space\number\day ,\number\year}

\newcommand{\proof}{\noindent{\bf Proof}~~}
\newtheorem{Theorem}{Theorem}

\newtheorem{Lemma}{Lemma}
\newcommand{\bl}{\begin{Lemma}}

\newcommand{\el}{\end{Lemma}}
\newcommand{\be}{\begin{equation}}
\newcommand{\ee}{\end{equation}}
\newcommand{\bd}{\begin{defi}}
\newcommand{\ed}{\end{defi}}

\newtheorem{pro}{Proposition}
\newcommand{\bp}{\begin{pro}}
\newcommand{\ep}{\end{pro}}
\newcommand{\bt}{\begin{Theorem}}
\newcommand{\et}{\end{Theorem}}
\newtheorem{cor}{Corollary}
\newcommand{\bc}{\begin{cor}}
\newcommand{\ec}{\end{cor}}

\def\sqr#1#2{{\vcenter{\vbox{\hrule height.#2pt
\hbox{\vrule width.#2pt height#1pt \kern#1pt
\vrule width.#2pt}\hrule height.#2pt}}}}
\def\square{\mathchoice\sqr45\sqr45\sqr{2.1}3\sqr{1.5}3}
\newcommand{\qed}{\square}

 at 10truept

\begin{document}
\setlength{\textheight}{7.7truein}    
\setcounter{page}{1} \centerline{\bf  D-homothetically fixed, weakly $(\kappa ,\mu )$-structures on contact metric  spaces } \centerline{\bf } \baselineskip=13pt
\vspace*{10pt} \centerline{\footnotesize {\bf 
Philippe Rukimbira}} \baselineskip=12pt \centerline{\footnotesize\it
Department of Mathematics and Statistics, Florida International University}
\baselineskip=10pt \centerline{\footnotesize\it Miami, Florida
33199, USA} \baselineskip=10pt \centerline{\footnotesize E-MAIL:
rukim@fiu.edu} \vspace*{0.225truein}

\vspace*{0.21truein} \abstract{\it } {Contact metric $(\kappa ,\mu )$-spaces are generalizations of Sasakian spaces. We introduce a weak $(\kappa ,\mu )$ condition as a generalization of the K-contact one and show that many of the known  results from generalized Sasakian geometry hold in the weaker generalized  K-contact geometry setting. In particular, we prove existence of K-contact and $(\kappa ,\mu =2)$-structures under some conditions on the Boeckx invariant.}


\vspace*{14pt}                  

\baselineskip=24pt
\section{Basic properties of contact metric structures}

A contact form on a $2n+1$-dimensional manifold $M$ is a one-form $\eta$ such that $\eta\wedge (d\alpha )^n$ is a volume form on $M$. Given a contact manifold $(M,\eta )$, there exist tensor fields $(\xi ,\phi , g)$, where $g$ is a Riemannian metric and $\xi$ is a unit vector field, called the Reeb field of $\eta$ and $\phi$ is an endomorphism of the tangent bundle of $M$ such that
\begin{itemize}
\item[(i)] $\eta (\xi )=1,~\phi^2 =-Id +\eta\otimes\xi ,~\phi\xi =0$
\item[(ii)] $d\eta =2g(.,\phi .)$
\end{itemize}
The data $(M,\eta ,\xi ,\phi ,g)$ is called a contact metric structure; see  (\cite{BLA}) for more details. 

Denoting by $\nabla$ the Levi-Civita connection of $g$, and by $$R(X,Y)Z=\nabla_X\nabla_YZ-\nabla_Y\nabla_XZ-\nabla_{[X,Y]}Z$$ its curvature tensor, a contact metric structure $(M,\eta ,\xi ,\phi ,g)$ is called Sasakian if the condition$$(\nabla_X\phi )Y=g(X,Y)\xi -\eta (Y)X$$ is satisfied for all tangent vectors  $X$ and $Y$.
A well known curvature characterization of the Sasakian condition is as follows:
\begin{pro}
{\it 
A contact metric structure $(M,\eta , \xi ,\phi ,g)$ is Sasakian if and only if $$R(X,Y)\xi =\eta (Y)X-\eta (X)Y$$ for all tangent vectors $X$ and $Y$.}\end{pro}

A condition weaker than the Sasakian one is the K-contact condition. A contact metric structure $(M,\eta ,\xi ,\phi ,g)$ is called K-contact if the tensor field $h=\frac{1}{2}L_\xi \phi$ vanishes identically. Here, $L_\xi \phi$ stands for the Lie derivative of $\phi$ in the direction of $\xi$. The above K-contact condition is known to be equivalent to the Reeb vector field $\xi$ being a $g$- infinitesimal isometry, or a Killing vector field. The tensor field $h$ is known to be symmetric and anticommutes with $\phi$.

An equally well known curvature characterization of K-contactness is as follows: 

\begin{pro}{\it A contact metric structure $(M,\eta ,\xi ,\phi ,g)$ is K-contact if and only if $$R(X,\xi )\xi =X-\eta (X)\xi$$ for all tangent vectors $X$.}\end{pro}

The notation ''$l$" is common for the tensor $$lX=R(X,\xi )\xi$$

\begin{pro}\label{prop33}

On a contact metric structure $(M,\eta ,\xi ,\phi ,g)$, the following identities hold:

\begin{equation}\nabla_\xi h=\phi -h^2\phi -\phi l \label{bl1}\end{equation}
\begin{equation}\phi l\phi -l=2(h^2+\phi^2) \label{bl2}\end{equation}
\begin{equation}L_\xi h=\nabla_\xi h+2\phi h+2\phi h^2\label{bl3}\end{equation}
\end{pro}

\proof  The first two identities appear in Blair's book (\cite{BLA}). We establish the third one.
$$\begin{array}{rcl}( L_\xi h) X&=&[\xi ,hX]-h[\xi ,X]\\&=&\nabla_\xi(hX)-\nabla_{hX}\xi -h(\nabla_\xi X-\nabla_X\xi )\\&=&(\nabla_\xi h)X+h\nabla_\xi X-[-\phi hX-\phi h^2X]-h\nabla_\xi X+h[-\phi X-\phi hX]\\
&=&(\nabla_\xi h)X+h\nabla_\xi X+\phi hX+\phi h^2X-h\nabla_\xi X-h\phi X+\phi h^2 X\\ &=&(\nabla_\xi h)X+2\phi hX+2\phi h^2X
\end{array}
$$

$\qed$

Given a contact metric structure $(M,\eta ,\xi , \phi , g)$, its $D_a$-homothetic deformation is a new contact metric structure $(M,\overline{\eta}, \overline{\xi }, \overline{\phi}, \overline{g})$ given by a real number $a>0$ and  $$\overline{\eta}=a\eta,~~\overline{\xi}=\frac{\xi}{a},~~\overline{\phi}=\phi$$
$$\overline{g}=ag+a(a-1)\eta\otimes\eta$$ 
D-homothetic deformations preserve the K-contact and Sasakian conditions.
\section{Weakly $(\kappa , \mu )$- spaces}

A direct calculation shows that under a $D_a$ homothetic deformation, the curvature tensor transforms as follows:
$$\begin{array}{rcl}a\overline{R}(X,Y)\overline{\xi}&=&R(X,Y)\xi -(a-1)[(\nabla_X\phi )Y-(\nabla_Y\phi )X+\eta (X)(Y+hY)\\&&-\eta (Y)(X+hX)]+(a-1)^2[\eta (Y)X-\eta (X)Y]
\end{array}
$$
Letting $Y=\xi$ and recalling $\nabla_\xi\phi =0$, we get:$$\begin{array}{rcl}a\overline{R}(X,\xi )\overline{\xi}&=&R(X,\xi )\xi-(a-1)[(\nabla_X\phi )\xi +\eta (X)\xi -(X+hX)\\&&+(a-1)^2[X-\eta (X)\xi ]
\end{array}
$$
On any contact metric manifold, the following identity holds: $$(\nabla_X\phi)\xi =-\phi\nabla_X\xi=-X+\eta (X)-hX.$$ 

Taking into account of this identity, we see that the curvature tensor deforms as follows:
$$a^2\overline{R}(X,\overline{\xi})\overline{\xi}=R(X,\xi )\xi +(a^2-1)(X-\eta (X)\xi )+2(a-1)hX$$

Equivalently, since $\xi =a\overline{\xi}$  and $h=a\overline{h}$, 
\begin{equation}\overline{R}(X,\overline{\xi})\overline{\xi}=\frac{1}{a^2}R(X,\xi )\xi +\frac{a^2-1}{a^2}(X-\overline{\eta} (X)\overline{\xi})+\frac{2a-2}{a}\overline{h}X\label{k6}\end{equation}

It follows from (\ref{k6}) that, 
under a $D_a$-homothetic deformation, the condition $R(X,\xi )\xi =0$ transforms into $$\overline{R}(X,\overline{\xi} )\overline{\xi} =\kappa (X-\overline{\eta} (X)\overline{\xi} )+\mu \overline{h}X$$ where $\kappa =\frac{a^2-1}{a^2}$ and $\mu =\frac{2a-2}{a}$.  

As a generalization of both $R(X,\xi )\xi =0$ and the K-contact condition, $R(X,\xi )\xi =X-\eta (X)\xi $, we consider $$R(X,\xi )\xi =\kappa (X-\eta(X)\xi ) +\mu hX.$$ We call this the weak $(\kappa , \mu )$ condition. The same generalization was refered to as Jacobi $(\kappa , \mu )$-contact s manifold in \cite{GHS}.  Let us point out also that  a strong $(\kappa ,\mu )$ condition $R(X,Y)\xi=\kappa (\eta(Y)X-\eta (X)Y)+\mu (\eta (Y)hX-\eta (X) hY)$ has been introduced in \cite{BKP}. Examples of weakly $(\kappa ,\mu )$ spaces which are not strongly $(\kappa ,\mu )$ are provided by the Darboux contact forms $\eta =\frac{1}{2}(dz-\sum y^idx^i)$ on ${\bf R}^{2n+1}$ with associated metric 
$$g=\frac{1}{4}\left(\begin{array}{ccr} \delta_{ij}+y^iy^j+\delta_{ij}z^2&\delta_{ij}z&-y^i\\\delta_{ij}z&\delta_{ij}&0\\-y^j&0&1\end{array}\right)
$$
(see \cite{BLA}).  

Other examples of weakly $(\kappa ,\mu )$-spaces have been found on normal bundles of totally geodesic Legendre submanifolds in Sasakian manifolds (see \cite{BAN}).

The two notions of $(\kappa ,\mu )$-spaces are D-homothetically invariant. It follows from identity (\ref{k6})  that,  if $(M,\eta ,\xi ,\phi ,g )$ is a (weak)  $(\kappa ,\mu )$ structure, then the $D_a$-homothetic deformation $( \overline{\eta}, \overline{\xi}, \phi ,\overline{g} )$  is a (weak) $(\overline{\kappa}, \overline{\mu})$-structure with:

 \begin{equation}\label{km} \overline{\kappa}=\frac{\kappa +a^2-1}{a^2},~~\overline{\mu}=\frac{\mu +2a-2}{a}\end{equation}
The tensor fields $\phi$ and $h$ on a weakly $(\kappa ,\mu )$-space are related by the identities in the following proposition.

\begin{pro}\label{prop4} On a weakly $(\kappa ,\mu )$-space $(M,\eta ,\xi , \phi , g)$, the following identities hold:
\begin{equation} h^2=(\kappa -1)\phi^2,~~\kappa \le 1\label{ka1}\end{equation}
\begin{equation}\nabla_\xi h=-\mu\phi h\label{ka2}\end{equation}
\begin{equation} L_\xi h=(2-\mu )\phi h+2(1-\kappa )\phi\label{eq6}\end{equation}
\end{pro}

\proof Starting with identity (\ref{bl2}), which is valid on any contact metric structure, one has, for any tangent vector $X$:
$$\begin{array}{rcl}2h^2X+2\phi^2X&=&\phi (\kappa \phi X+\mu h\phi X)-(\kappa (-\phi^2X)+\mu hX)\\&=&\kappa\phi^2X+\mu\phi h\phi X+\kappa\phi^2X-\mu hX\\&=&\kappa\phi^2X-\mu h\phi^2X+\kappa\phi^2X-\mu hX\\&=&2\kappa\phi^2X\end{array}
$$
Hence, grouping terms
$$2h^2X=(2\kappa -2)\phi^2X$$ So $$h^2=(\kappa -1)\phi^2$$ 

But, since $h$ is symmetric, $h^2$ must be a non-negative operator, hence $\kappa \le 1$, proving (\ref{ka1}).

 From identity (\ref{bl1}) combined with $lX=\kappa (X-\eta (X)\xi )+\mu hX$, we see that 
$$\begin{array}{rcl}(\nabla_\xi h)X&=&\phi X-h^2\phi X-\phi (\kappa (X-\eta (X)\xi )+\mu hX)\\&=&\phi X-(\kappa -1)\phi^3X-\kappa\phi X-\mu\phi hX\\&=&(1-\kappa +(\kappa -1))\phi X-\mu\phi hX\\&=&-\mu \phi hX
\end{array}
$$  proving (\ref{ka2}).

Next, combining identities (\ref{bl3}), (\ref{ka2}) and (\ref{ka1}), one has:
$$\begin{array}{rcl}L_\xi h&=&\nabla_\xi h+2\phi h+2\phi h^2\\&=&-\mu \phi h+2\phi h+2\phi h^2\\
&=&-\mu \phi h+2\phi h+2\phi (\kappa -1)\phi^2\\&=&-\mu \phi h+2\phi h-2(\kappa -1)\phi\\&=&(2-\mu )\phi h+2(1-\kappa )\phi\end{array}
$$ proving (\ref{eq6}).
$\qed$

Tangent bundle's structure on a $(\kappa ,\mu )$-space is described by the following theorem:
\begin{Theorem} Let $(M^{2n+1},\eta , \xi ,\phi , g)$ be a weakly $(\kappa ,\mu )$, contact metric manifold. Then $\kappa \le 1$. If $\kappa =1$, then the structure is K-contact. If $\kappa <1$, then the tangent bundle $TM$ decomposes into three mutually orthogonal distributions $D(0)$, $D(\lambda )$ and $D(-\lambda )$, the eigenbundles determined by tensor $h$'s eigenspaces, where $\lambda =\sqrt{1-\kappa }$.
\end{Theorem} 

\proof  Clearly, $\kappa =1$ is exactly the K-contact condition.

Suppose $\kappa <1$. Since $h\xi =0$ and $h$ is symmetric, it follows from identity (\ref{ka1}), Proposition \ref{prop4}, ($h^2=(\kappa -1)\phi^2$),  that the restriction $h_{|D}$ of $h$ to the contact subbundle $D$ has eigenvalues $\lambda =\sqrt{1-k}$ and $-\lambda$.    By $D(\lambda ),~D(-\lambda )~and~D(0)$, we denote the corresponding eigendistributions. If $X\in D(\lambda )$, then $h\phi X=-\phi hX=-\lambda\phi X$. Thus $\phi X\in D(-\lambda )$ which shows that the three distributions above are mutually orthogonal.   $\qed$

To shade some light on  the difference between weak $(\kappa ,\mu )$ and strong $(\kappa ,\mu )$-spaces, we propose a weak, semi-symmetry condition. We say that a contact metric space $(M,\eta ,\xi ,\phi , g)$ is weakly semi-symmetric if $R(X,\xi )R=0$ for all tangent vectors $X$ where $R$ is the curvature operator.

We will prove the following:

\begin{Theorem}\label{Theo2} Let $(M,\eta ,\xi ,\phi , g)$ be a weakly semi-symmetric,  contact metric weakly $(\kappa ,0)$-space. Then $(M,\eta ,\xi ,\phi ,g )$ is a strongly $(\kappa ,0)$-space.
\end{Theorem} 
\proof

The weakly semi- symmetric condition means that $(R(X,\xi )R)(Y,\xi )\xi =0$ holds for any tangent vectors $X$ and $Y$. Extending $Y$ into a local vector field, we have:
$$\begin{array}{rcl} 0&=&R(X,\xi )R(Y,\xi )\xi-R(R(X,\xi )Y,\xi )\xi -R(Y,R(X,\xi )\xi )\xi -\\&&R(Y,\xi )R(X,\xi )\xi\\

&=&R(X,\xi )(\kappa (Y-\eta (Y)\xi ))-\kappa (R(X,\xi )Y-\eta (R(X,\xi )Y)\xi )-\\&&R(Y,\kappa (X-\eta (X)\xi ))\xi -R(Y,\xi )(\kappa (X-\eta (X)\xi ))\\

&=&\kappa R(X,\xi )Y-\kappa \eta (Y)R(X,\xi )\xi-\kappa R(X,\xi )Y+\kappa\eta (R(X,\xi)Y)\xi -\\&&\kappa R(Y,X)\xi +\kappa\eta (X)R(Y,\xi )\xi -\kappa R(Y,\xi )X+\kappa\eta (X)R(Y,\xi )\xi\\

&=&-\kappa^2\eta (Y)X-\kappa g(R(X,\xi )\xi ,Y)\xi -\kappa R(Y,X)\xi +\kappa^2\eta (X)Y-\\&&\kappa R(Y,\xi )X+\kappa^2\eta (X)Y-\kappa^2\eta (X)\eta (Y)\xi\\&=&-\kappa^2\eta (Y)X-\kappa g(\kappa (X-\eta (X)\xi ),Y)\xi -\kappa R(Y,X)\xi +2\kappa^2\eta (X)Y-\\&&\kappa R(Y,\xi )X-\kappa^2\eta (X)\eta (Y)\xi\\0&=&-\kappa^2\eta (Y)X-\kappa^2 g(X,Y)\xi -\kappa R(Y,X)\xi +2\kappa^2\eta (X)Y-\kappa R(Y,\xi )X\end{array}
$$Equation \begin{equation}\label{w1} -\kappa^2\eta (Y)X-\kappa^2 g(X,Y)\xi -\kappa R(Y,X)\xi +2\kappa^2\eta (X)Y-\kappa R(Y,\xi )X=0\end{equation} is valid for any $X$ and $Y$. Exchanging  $X$ and $Y$ leads to
\begin{equation}\label{w2}-\kappa^2\eta (X)Y-\kappa^2g(Y,X)\xi -\kappa R(X,Y)\xi +2\kappa^2\eta (Y)X-\kappa R(X,\xi )Y=0\end{equation}
Substracting equation (\ref{w1}) from equation (\ref{w2}), we obtain:\begin{equation}3\kappa^2 (\eta (Y)X-\eta (X)Y)+\kappa (R(Y,X)\xi -R(X,Y)\xi  )+\kappa (R(Y,\xi )X-R(X,\xi )Y)=0\label{w3}\end{equation}

By the first Bianchi Identity, $R(Y,\xi )X-R(X,\xi )Y=R(Y,X)\xi$ holds. Incorporating this identity into (\ref{w3}), we obtain the following:
\begin{equation} 3 \kappa^2(\eta (Y)X-\eta (X)Y)+3\kappa R(Y,X)\xi =0\end{equation} which implies the strong $(\kappa ,0)$ condition $$R(X,Y)\xi =\kappa (\eta (Y)X-\eta (X)Y)$$
$\qed$

Theorem \ref{Theo2} applies to the case $\kappa =1$ and has the following interesting corollary:
\begin{cor} A weakly semi-symmetric K-contact manifold is Sasakian.
\end{cor}

\proof  In the $\kappa =1$ case, the strong $(\kappa ,\mu )$ condition is exactly the Sasakian condition.   $\qed$

Weakly $(\kappa ,\mu )$-structures with $\kappa =1$, (the $K$-contact ones), are D-homothetically fixed. A non-K-contact weakly $(k,\mu )$ structure cannot be D-homothetically deformed into a K-contact one. Weakly $(k,\mu )$-structures with $\mu =2$ are also D-homothetically fixed. As a consequence, a  weakly $(k, \mu )$ structure with $\mu \neq 2$ cannot be deformed into one with $\mu =2$ neither.

Existence of these homothetically fixed $(\kappa ,\mu )$ structures  depends on  an invariant that was first introduced by Boeckx for strongly $(\kappa ,\mu )$ structures in (\cite{BOE}).

\section{The Boeckx invariant}
The Boeckx invariant, $I_M$ of a weakly contact $(\kappa ,\mu )$ space is defined by $$I_M=\frac{1-\frac{\mu }{2}}{\sqrt{1-k}}=\frac{1-\frac{\mu}{2}}{\lambda}.$$   $I_M$ is a D-homothetic invariant. Any two D-homothetically related weakly $(\kappa ,\mu )$-structures have the same Boeckx invariant.

The following lemma is crucial in proving existence of $D$-homothetically fixed $(\kappa ,\mu )$-structures.
\bl\label{lemma1} Let $(M,\alpha , \xi ,\phi , g )$ be a non-K-contact, weakly $(\kappa ,\mu )$ space.\begin{itemize}
\item[(i)] If $I_M>1$, then $2-\mu -\sqrt{1-k}>0$ and $2-\mu +\sqrt{1-k}>0$.
\item[(ii)] If $I_M<-1$, then $2-\mu +\sqrt{1-k}<0$ and $2-\mu -\sqrt{1-k}<0.$
\item[(iii)] $|I_M|<1$ if and only if $0<2\lambda +2-\mu$ and $0<2\lambda +\mu -2$
\end{itemize}
\el
 
\proof  (i).    Suppose $I_M>1$. Then $1-\frac{\mu}{2}>\sqrt{1-\kappa}$ and $\mu <2$.
$$\begin{array}{rcr}1-\frac{\mu}{2} >\sqrt{1-\kappa}&\Rightarrow&2-\mu >2\sqrt{1-k}\\&\Rightarrow& 2-\mu >\sqrt{1-k} ~and~ 2-\mu >-\sqrt{1-k}\\&\Rightarrow&2-\mu-\sqrt{1-k}>0~and~2-\mu +\sqrt{1-k}>0\end{array}
$$

(ii). Suppose $I_M<-1$. Then $1-\frac{\mu }{2}<-\sqrt{1-k}$ and $\mu >2$.
$$\begin{array}{rcr} 1-\frac{\mu}{2}<-\sqrt{1-k}&\Rightarrow&2-\mu <-2\sqrt{1-k} ~and ~2-\mu <\sqrt{1-k}\\&\Rightarrow&2-\mu <-\sqrt{1-k}~and ~2-\mu <\sqrt{1-k}\\&\Rightarrow&2-\mu +\sqrt{1-k}<0~and~2-\mu -\sqrt{1-k}<0\\&&\qed\end{array}
$$ 

(iii).  $|I_M|<1$ if and only if $-1<\frac{1-\frac{\mu}{2}}{\lambda}<1$. Equivalently
$$-1<\frac{2-\mu}{2\lambda}<1~~and~~-1<\frac{\mu -2}{2\lambda}<1$$

Thus $$-2\lambda <2-\mu <2\lambda~~and~~-2\lambda <\mu -2<2\lambda$$ Or, $$0<2\lambda +2-\mu ~and~~0<2\lambda +\mu -2$$
$\qed$

\section{D-homothetically fixed structures on weakly $(\kappa ,\mu )$-spaces}
\subsection{K-contact structures on weakly $(\kappa ,\mu )$ spaces}
We have pointed out that D-homothetic deformations of  non K-contact weakly $(\kappa ,\mu )$ structures remain non K-contact. However, on $(\kappa ,\mu )$-spaces with large Boeckx invariant, K-contact structures coexist with $(\kappa ,\mu )$ structures. 

\bt \label{theo2}Let $(M,\eta , \xi , \phi  ,g )$ be a non-K-contact, weakly $(k,\mu )$-space whose Boeckx invariant $I_M$ satisfies $|I_M|>1$. Then, $M$ admits a K-contact structure $(M,\eta , \xi, \overline{\phi}, \overline{g})$ compatible with the contact form $\eta$.\et
\proof  We define tensor fields $\overline{\phi}$ and $\overline{g}$ by   \begin{equation}\label{def1}\overline{\phi}=\frac{\epsilon}{(1-k)\sqrt{(2-\mu )^2-4(1-k)}}(L_\xi h\circ h)\end{equation}

$$\overline{g}=-\frac{1}{2}d\eta (., \overline{\phi}.)+\eta\otimes\eta$$ where $$\epsilon=\left\{\begin{array}{ll}+1&if~I_M>0\\-1&if~I_M<0\end{array}\right.$$

From the formula $h^2=-(1-k)\phi^2$  and $L_\xi h=(2-\mu )\phi h+2(1-k)\phi$ in Proposition \ref{prop4}, we obtain $$\begin{array}{rcl}
(L_\xi h\circ h)^2&=&(2-\mu )^2(1-k)^2\phi^2-4(1-k)^2\phi^2h^2\\
&=&(1-k)^2((2-\mu )^2-4(1-k))\phi^2
\end{array}
$$  
That is:  $$(L_\xi h\circ h)^2=\lambda^4\alpha (-Id+\eta\otimes \xi )$$ where $\lambda =\sqrt{1-k}$ and $\alpha =(2-\mu )^2-4(1-k).$ One sees that, if $\alpha >0$, then  $\overline{\phi}=\frac{\epsilon}{\lambda^2\sqrt\alpha}(L_\xi h\circ h)$ defines an almost complex structure on the contact subbundle. Notice also that $\alpha >0$ is equivalent to $|I_M|>1$.

We will show that $\overline{\phi}$ is $\xi$ invariant. For that, it suffices to show that the Lie derivative of $L_\xi h\circ h$ vanishes in the $\xi$ direction.
$$\begin{array}{rcl}L_\xi (L_\xi h\circ h)&=&L_\xi ((2-\mu )(1-k)\phi +2(1-k)\phi h)\\&=&2(2-\mu )(1-k )h+4(1-k)h^2+2(1-k)[(2-\mu )\phi^2h+\\&&2(1-k)\phi^2]\\
&=&2(2-\mu )(1-k)h+4(1-k)h^2+2(1-k)(2-\mu )\phi^2h+\\&&4(1-k)^2\phi^2\\
&=&-4(1-k)^2\phi^2+4(1-k)^2\phi^2 =0
\end{array}
$$
Next, we will show that $\overline{g}=-\frac{1}{2}d\eta (., \overline{\phi}. )+\eta\otimes\eta$ is an adapted Riemannian  metric for the structure tensors $(\eta , \overline{\xi}, \overline{\phi})$. That is $\overline{g}$ is a bilinear, symmetric, positive definite tensor with $$d\eta =2 \overline {g}(., \overline{\phi}).$$
 From the definition of $\overline{g}$, we have, for arbitrary tangent vectors $X$ and $Y$:

$$\begin{array}{rcl} \overline{g}(X,Y)&=& -\frac{1}{2}d\eta (X\overline{\phi}Y)+\eta (X)\eta (Y)\\
&=&-\frac{1}{2(1-\kappa )\sqrt{4(1-\kappa )-(2-\mu )^2}}d\eta (X,(1-\kappa )(2-\mu )\phi +\\&&(1-\kappa )2\phi h)Y)+\eta (X)\eta (Y)\\&=&
\frac{2-\mu }{\sqrt{4(1-\kappa )-(2-\mu )^2}}g(X,Y)+\frac{2}{\sqrt{4(1-\kappa )-(2-\mu )^2}}g(X,hY)+\\&&(1-\frac{1}{\sqrt{4(1-\kappa )-(2-\mu )^2}})\eta(X)\eta(Y)\\&=&\frac{2-\mu }{\sqrt{4(1-\kappa )-(2-\mu )^2}}g(Y,X)+\frac{2}{\sqrt{4(1-\kappa )-(2-\mu )^2}}g(Y,hX)+\\&&(1-\frac{1}{\sqrt{4(1-\kappa )-(2-\mu )^2}})\eta(Y)\eta(X\\&=&\overline{g}(Y,X)
\end{array}
$$ 
proving symmetry of $\overline{g}$. We used $h$'s symmetry in the step before the last.

For $\overline{g}$'s positive definitness, first observe that $\overline{g}(\xi ,\xi )=1>0$. Then for any non-zero tangent vector $X$ in the contact bundle $D$, using the definition of $\overline{\phi }$ in (\ref{def1}), the formula for $L_\xi h$ from identity (\ref{eq6}) in Proposition \ref{prop4}, we have:
$$\begin{array}{rcl}\overline{g}(X,X)&=&-\frac{1}{2}d\eta (X,\overline{\phi }X)\\&=&-\frac{\epsilon (2-\mu )}{2\sqrt{(2-\mu )^2-4(1-\kappa )}}d\eta (X,\phi X)-\frac{\epsilon}{\sqrt{(2-\mu )^2-(4(1-\kappa )}}d\eta (X,\phi hX)\\
&=&\frac{\epsilon (2-\mu )}{\sqrt{(2-\mu )^2-4(1-\kappa )}}g(X,X)+\frac{2\epsilon}{\sqrt{(2-\mu )^2-4(1-\kappa )}}g(X,hX)\\&=&
\frac{\epsilon}{\sqrt{(2-\mu )^2-4(1-\kappa )}}((2-\mu )g(X,X)+2g(X,hX))
\end{array}
$$
If $X\in D(\lambda )$, then $$\overline{g}((X,X)=\frac{\epsilon}{\sqrt{(2-\mu )^2-4(1-\kappa )}}((2-\mu )+2\sqrt{1-\kappa })g(X,X)).$$ By Lemma \ref{lemma1}, (i), (ii), the inequality $$\epsilon ((2-\mu )-2\sqrt{1-\kappa}))>0$$ holds when $|I_M|>1$.  Therefore $\overline{g}(X,X)>0$. 

In the same way,  if $X\in D(-\lambda )$, then $$\overline{g}(X,X)=\frac{\epsilon}{\sqrt{(2-\mu )^2-4(1-\kappa )}}((2-\mu )-2\sqrt{1-\kappa })g(X,X))$$
which is also $>0$ by Lemma \ref{lemma1}, (i) and (ii). This concludes the proof of $\overline{g}$'s positivity.

We easily verify that $\overline{g}$ is an adapted metric.

$$\begin{array}{rcl}2\overline{g}(X,\overline{\phi}Y)&=&-d\eta (X,\overline{\phi}^2Y)\\&=&-d\eta (X,-Y+\eta (Y)\xi )\\&=&d\eta (X,Y).\end{array}
$$  $\qed$

\noindent{\bf Remark:} As a consequence of Theorem \ref{theo2}, contact forms on compact, weakly $(\kappa ,\mu )$-spaces with $|I_M|>1$ admit associated K-contact structures,  hence verify Weinstein's Conjecture about the existence of closed Reeb orbits.( see\cite{ RUK}).

On weakly $(\kappa ,\mu )$-spaces with small Boeckx invariant, it turns out that $(\kappa ,2)$ structures coexist with $(\kappa ,\mu\ne 2)$ structures. This will be established in the next subsection.

\subsection{Contact metric weakly  $(\kappa , 2 )$-spaces}

Given a non-K-contact, weakly $(\kappa ,\mu )$-space $(M,\eta ,\xi ,\phi ,g)$, we define the D-homothetic invariant tensor field $$\tilde{\phi}=\frac{1}{\sqrt{1-k}}h$$

\begin{Lemma}\label{lem2} Denoting by $$\tilde{h}=\frac{1}{2}L_\xi\tilde{\phi }=\frac{1}{2\sqrt{1-k}}L_\xi h,$$ the following identities are satisfied:
\begin{equation}\tilde{h}=\frac{1}{2\sqrt{1-k}}((2-\mu )\phi h+2(1-k)\phi )\label{tilde1}\end{equation}
\begin{equation}\tilde{h}^2=((1-k)-(1-\frac{\mu}{2})^2)\phi^2 \label{tilde2}\end{equation}
\end{Lemma}

\proof   From the third identity in Proposition \ref{prop33}, combined with identity (\ref{eq6}), Proposition \ref{prop4}, we get       $$2(\sqrt{1-k})\tilde{h}=L_\xi h=(2-\mu )\phi h+2(1-\kappa )\phi $$  

So

$$\tilde{h}=\frac{1}{2\sqrt{1-k}}(2-\mu )\phi h+2(1-\kappa )\phi$$
 which is 
(\ref{tilde1}).

The proof of (\ref{tilde2}) is a straightforward calculation.  $\qed$

\noindent {\bf Remark}: If $|I_M|<1$, then $1-k-(1-\frac{\mu}{2})^2>0$. Therefore, identity (\ref{tilde2}) suggests that  $\tilde{h}$ can be used to define a complex structure on the contact subbundle.\vskip 12pt
Define  the tensor field $\phi_1$ by:$$\phi_1=\frac{1}{\sqrt{1-k-(1-\frac{\mu}{2})^2}}\tilde{h}=\frac{1}{\sqrt{1-k-(1-\frac{\mu}{2})^2}}\frac{1}{2\sqrt{1-k}}((2-\mu)\phi h+2(1-k)\phi )$$

\begin{pro} \label{prop5}The tensor field $\phi_1$ satisfies $$\phi_1^2=-I+\eta\otimes \xi$$ 
\begin{equation} h_1=\frac{1}{2}L_\xi\phi_1=(\sqrt{1-I_M^2})h\label{phi11}\end{equation}\end{pro}

\proof The identity $\phi_1^2=\phi^2=-I+\eta\otimes\xi$ follows from  Lemma \ref{lem2}, (\ref{tilde2}). 

As for identity (\ref{phi11}), we proceed as follows: 
$$\begin{array}{rcl}h_1=\frac{1}{2}(L_\xi\phi_1)&=&\frac{1}{4\sqrt{(1-\kappa )(1-\kappa-(1-\frac{\mu}{2})^2}}
L_\xi ((2-\mu )\phi h+2(1-\kappa )\phi\\&=&\frac{1}{4 \sqrt{(1-\kappa )(1-\kappa -(1-\frac{\mu}{2})^2}}[(2-\mu)((L_\xi\phi )h+\phi L_\xi h)+\\&&2(1-\kappa )L_\xi \phi]\\&=& \frac{1}{4 \sqrt{(1-\kappa )(1-\kappa -(1-\frac{\mu}{2})^2}}((2-\mu )(2h^2+\phi((\mu -2)h\phi+\\&&2(1-\kappa )\phi ) +4(1-\kappa )h)\\
&=&\frac{1}{4 \sqrt{(1-\kappa )(1-\kappa -(1-\frac{\mu}{2})^2}}[2(2-\mu )h^2-(2-\mu )^2\phi h\phi +\\&&2(2-\mu )(1-\kappa )\phi^2+4(1-\kappa )h]\\&=&\frac{1}{4 \sqrt{(1-\kappa )(1-\kappa -(1-\frac{\mu}{2})^2}}[2(2-\mu )(\kappa -1)\phi^2-(2-\mu )^2h+\\&&2(2-\mu )(1-\kappa )\phi^2 +4(1-\kappa )h]\\&=&  \frac{1}{4 \sqrt{(1-\kappa )(1-\kappa -(1-\frac{\mu}{2})^2}}[4(1-\kappa )-(2-\mu )^2]h\\&=&\sqrt{\frac{4(1-\kappa )-(2-\mu )^2}{4(1-\kappa )}}h=(\sqrt{1-I_M^2})h    \end{array}
$$
$\qed$

As pointed out earlier, 
when a D-homothetic deformation is applied to a weakly $(\kappa , \mu )$ structure with $\mu =2$, the $\mu$  value remains the same as is seen from one of formulas (\ref{km}):
$$\overline{\mu}=\frac{\mu +2a-2}{a}$$
 As a consequence,  weakly $(\kappa ,2 )$ structures cannot be obtained through D-homothetic deformations. In the case $|I_M|<1$, we prove  the following theorem:
\bt  \label{theo4} Let $(M, \eta , \xi , \phi , g)$ be a non-K-contact,  weakly $(\kappa , \mu )$-space with Boeckx invariant $I_M$ satisfying  $ |I_M|<1$. Then, there is a weakly $(\kappa_1 ,\mu_1 )$ structure $(M,\eta , \xi , \phi_1, g_1 )$ where $\mu_1=2$ and $\kappa_1=\kappa +(1-\frac{\mu}{2})^2$.\et

\proof  Define $g_1$ by  $$g_1(X,Y)=-\frac{1}{2}d\eta (X,\phi_1 Y)+\eta (X)\eta (Y).$$  We will show that $g_1$ is a Riemannian metric adapted to $\phi_1$ and $\eta$, i.e. $$d\eta =2g_1(.,\phi_1 )$$
For any tangent vectors $X$ and $Y,$ 
$$\begin{array}{rcl}g_1(X,Y)&=&-\frac{1}{2}\frac{1}{\sqrt{1-\kappa -(1-\frac{\mu}{2})^2}}d\eta (X,\tilde{h}Y)+\eta (X)\eta (Y)\\&=&-\frac{1}{2\sqrt{1-\kappa -(1-\frac{\mu}{2})^2}}d\eta (X,(2-\mu)\phi hY)+d\eta (X,2(1-\kappa )\phi Y))+\\&&\eta (X)\eta (Y)\\&=&-\frac{1}{2\sqrt{1\kappa  -(1-\frac{\mu}{2})^2}}(2g(X,(2-\mu )\phi^2hY)+2d\eta (X, 2(1-\kappa )\phi^2Y))+\\&&\eta (X)\eta (Y)\\&=& -\frac{1}{2\sqrt{1-\kappa -(1-\frac{\mu}{2})^2}}d\eta (X,(2-\mu)\phi hY)+d\eta (X,2(1-\kappa )\phi Y))+\\&&\eta (X)\eta (Y))\\&=&-\frac{1}{2\sqrt{1-\kappa -(1-\frac{\mu}{2})^2}}(2g((2-\mu )\phi^2hX,Y)+2d\eta (2(1-\kappa )\phi^2X,Y))+\\&&\eta (Y)\eta (X)\\&=& -\frac{1}{2}\frac{1}{\sqrt{1-\kappa -(1-\frac{\mu}{2})^2}}d\eta (Y,\tilde{h}X)+\eta (Y)\eta (X)\\&=&  g_1(Y,X)\end{array}
$$ proving that $g_1$ is a symmetric tensor.
To prove positivity of $g_1$, first observe that $g_1(\xi ,\xi )=1>0$. Next, for any $X$ in the contact distribution,
$$\begin{array}{rcl}g_1(X,X)&=&-\frac{1}{2}d\eta (X,\frac{1}{2\sqrt{1-\kappa}\sqrt{1-\kappa -(1-\frac{\mu}{2})^2}}((2-\mu )\phi hX+2(1-\kappa ) \phi X)\\&=&-\frac{(2-\mu )}{4\sqrt{1-\kappa }\sqrt{1-\kappa -(1-\frac{\mu}{2})^2}}d\eta (X,\phi hX)-\frac{(1-\kappa )}{2\sqrt{1-\kappa}\sqrt{1-\kappa -(1-\frac{\mu}{2})^2}}d\eta (X,\phi X)\\&=&\frac{1}{2\sqrt{1-\kappa }\sqrt{1-\kappa -(1-\frac{\mu}{2})^2}}((2-\mu )g(X,hX)+2(1-\kappa )g(X,X))\end{array}
$$
\begin{equation}\label{fr1}g_1(X,X)=\frac{1}{2\sqrt{1-\kappa }\sqrt{1-\kappa -(1-\frac{\mu}{2})^2}}((2-\mu )g(X,hX)+2(1-\kappa )g(X,X))\end{equation}

If $X\in D(\lambda )$, then (\ref{fr1}) becomes $$g_1(X,X)=\frac{1}{\sqrt{1-\kappa -(1-\frac{\mu}{2})^2}}(2\sqrt{1-\kappa }+(2-\mu ))g(X,X))>0 $$
The last inequality follows from Lemma \ref{lemma1}, (iii).
If $X\in D(-\lambda )$, then (\ref{fr1}) becomes $$g_1(X,X)=\frac{1}{\sqrt{1-\kappa-(1-\frac{\mu}{2})^2}} (2\sqrt{1-\kappa} -(2-\mu ))g(X,X)>0$$ also following from Lemma \ref{lemma1}, (iii).

We now prove that $g_1$ is an adapted metric.  Directly from the definition of $g_1$, 
$$\begin{array}{rcl}2g_1(X, \phi_1Y)&=&-d\eta (X,\phi^2_1Y)\\&=&d\eta (X,Y)
\end{array}
$$  

$\qed$

Finally, we show that the structure $(M,\eta ,\xi ,\phi_1, g_1)$ is a weakly $(\kappa_1 ,2 )$-structure. By Proposition \ref{prop5}, (\ref{phi11}), the positive eigenvalue of $h_1$ is $$\lambda_1=\sqrt{1-I_M^2}\lambda =\sqrt{(1-\kappa )(1-I_M^2)}=\sqrt{1-\kappa -(1-\frac{\mu}{2})^2}.$$ Since $(\eta ,\xi , \phi_1, g_1)$ is a contact metric structure, identity (\ref{bl1}), Proposition \ref{prop33} holds.
$$\overline{\nabla}_\xi h_1=\phi_1-\phi_1l_1-\phi_1h_1^2.$$
For any tangent vector field $X$, one has
$$\phi_1X-\phi_1l_1X-\phi_1h^2X=(\overline{\nabla}_\xi h_1)X$$
$$\begin{array}{rcl}\phi_1X-\phi_1l_1X-\lambda_1^2\phi_1X&=&\overline{\nabla}_\xi (h_1X)-h_1\overline{\nabla}_\xi X\\&=&\overline{\nabla}_{h_1X}\xi +[\xi , h_1X]-h_1(\overline{\nabla}_X\xi +[\xi ,X])\\&=&-\phi_1h_1X-\phi_1h_1^2X+(L_\xi h_1)X+h_1[\xi ,X]\\&&-h_1(-\phi_1X-\phi_1h_1X+[\xi ,X])\\ \phi_1X-\phi_1l_1X-\lambda_1^2\phi_1X&=&-2\phi_1h_1X-2\lambda_1^2\phi_1X+(L_\xi h_1)X 
\end{array}
$$
Applying $\phi_1$ on both sides of the above identity, one has
$$\phi_1^2X+l_1X-\lambda_1^2\phi_1^2X=2h_1X-2\lambda_1^2\phi_1^2X+\phi_1(L_\xi h_1)X$$
Solving for the tensor field $l_1$ gives
\begin{equation}\label{m2}l_1X=2h_1X-(1+\lambda_1^2)\phi_1^2X+(\phi_1L_\xi h_1)X\end{equation}

From Proposition \ref{prop5}, we know $L_\xi h_1=\sqrt{1-I_M^2}L_\xi h$ and $L_\xi h=(\mu -2)h\phi +2(1-\kappa )\phi$ from Proposition \ref{prop4}. Also $\phi_1=\frac{1}{2\sqrt{(1-\kappa )(1-\kappa -(1-\frac{\mu}{2})^2}}(L_\xi h)$. A direct calculation shows that
$$\begin{array}{rcl} \phi_1L_\xi h_1&=&\frac{\sqrt{1-I_M^2}}{2\sqrt{(1-\kappa )(1-\kappa -(1-\frac{\mu}{2})^2}}(L_\xi h)^2\\
&=&\frac{\sqrt{1-I_M^2}}{2\sqrt{(1-\kappa )(1-\kappa -(1-\frac{\mu}{2})^2}}(1-\kappa )[4(1-\kappa )-(\mu -2)^2]\phi^2\\&=&\frac{1}{2}[4(1-\kappa )-(\mu -2)^2]\phi^2\end{array}
$$
Reporting this in identity (\ref{m2}), we get:
$$\begin{array}{rcl}l_1X&=&2h_1X-(1+\lambda_1^2)\phi_1^2X+\frac{1}{2}[4(1-\kappa )-(2-\mu )^2]\phi^2X\\&=&2h_1X-(1+1-\kappa -(1-\frac{\mu}{2})^2)(-X+\eta (X)\xi +2(1-\kappa )-\\&&\frac{(2-\mu )^2}{2}(-X+\eta (X)\xi\\&=&2h_1X+(\kappa +\frac{(2-\mu )^2}{4})(X-\eta (X)\xi)
\end{array}
$$ Which is the $(\kappa_1, \mu_1)$ condition with $\mu_1=2$ and $\kappa_1=\kappa +\frac{(2-\mu)^2}{4}$.
$\qed$

\end{document}